\documentclass[twoside,twocolumn,10pt]{article}

\usepackage{wscg}           
\usepackage{algorithm}
\usepackage{algpseudocode}
\usepackage{hyperref}       
\usepackage{booktabs}       
\usepackage{amsfonts}       
\usepackage{longtable}
\usepackage{amsmath}
\usepackage{amssymb} 

\RequirePackage{ifpdf}
\ifpdf
 \RequirePackage[pdftex]{graphicx}
 \RequirePackage[pdftex]{color}
\else
 \RequirePackage[dvips,draft]{graphicx}
 \RequirePackage[dvips]{color}
\fi

\usepackage{nopageno}       

\title{Hybrid Quantum Branch-and-Bound Method for Quadratic Unconstrained Binary Optimization}

\author{
\parbox{0.25\textwidth}{\centering
Zedong Peng\\[1mm]
Purdue University\\
610 Purdue Mall\\
47907, West Lafayette, Indiana, USA\\[1mm]
peng372@purdue.edu
}
\hspace{0.05\textwidth}
\parbox{0.25\textwidth}{\centering
Daniel de Roux\\[1mm]
Carnegie Mellon University\thanks{blahag}\\
5000 Forbes Ave\\
15213, Pittsburgh, Pennsylvania, USA\\[1mm]
dderoux@andrew.cmu.edu\\[1mm]
*now at Google.
}
\hspace{0.05\textwidth}
\parbox{0.25\textwidth}{\centering
David E. Bernal Neira\\[1mm]
Purdue University\\
610 Purdue Mall\\
47907, West Lafayette, Indiana, USA\\[1mm]
dbernaln@purdue.edu
}
}

\usepackage{url}
\makeatletter
\def\Uslash{\mathbin{\mathchar`\/}\@ifnextchar{/}{\kern-.15em}{}}
\g@addto@macro\UrlSpecials{\do \/ {\Uslash}}
\def\Ucolon{\mathbin{\mathchar`:}\@ifnextchar{/}{\kern-.1em}{}}
\g@addto@macro\UrlSpecials{\do : {\Ucolon}}
\makeatother

\begin{document}

\twocolumn[{\csname @twocolumnfalse\endcsname

\maketitle  

\begin{abstract}
\noindent
Quantum algorithms have shown promise in solving Quadratic Unconstrained Binary Optimization (QUBO) problems, benefiting from their connection to the transverse field Ising model. Various Ising solvers, both classical and quantum, have emerged to tackle such problems efficiently but lack global optimality guarantees and often suffer from hardware limitations such as limited qubit availability. In this work, we propose a hybrid branch-and-bound (B\&B) framework that integrates Ising solvers as heuristics within a classical B\&B algorithm. Unlike prior theoretical studies, our work presents a practical implementation, available as open-source on GitHub. We explore when and where to apply Ising solvers in the search tree and introduce a custom branching rule optimized QUBO embedding. Our method is evaluated on hundreds of QUBO instances from QUBOLib.jl using Gurobi and the D-Wave quantum annealer. Our results show up to 11\% less solution time and 17\% fewer nodes compared to default Gurobi, an off-the-shelf commercial optimization solver. These findings demonstrate the value of hybrid quantum-classical strategies for enhancing exact optimization.

\end{abstract}

\subsection*{Keywords}
Quadratic Unconstrained Binary Optimization, Branch-and-Cut, Quantum Optimization.
\vspace*{1.0\baselineskip}
}]


\newcommand{\db}[1]{\textcolor{red}{#1}}


\section{Introduction}

The quadratic unconstrained binary optimization problem (QUBO) is a central family of optimization programs where one seeks to minimize a quadratic function over binary variables.  Formally, a QUBO is defined as
\begin{equation}
\label{qubo}
\begin{aligned}
& \text{Minimize}   && \mathbf{x}^\top \mathbf{Q} \mathbf{x} \\
& \text{subject to} && x_i \in \{0, 1\}, && \forall i \in N,
\end{aligned}
\end{equation}
where $\mathbf{x} \in \{0,1\}^n$, $N = \{1,\dots,n\}$, and $\mathbf{Q} \in \mathbb{R}^{n \times n}$ is symmetric.  
The objective function includes both linear terms $Q_{ii}x_i$ and bilinear terms $Q_{ij}x_ix_j$.

This family of optimization problems arises in various applications, such as MAXCUT, graph partitioning, protein folding, machine learning pipelines, and many others~\cite{QUBOHandbook}. QUBO is a special case of mixed-integer programming (MIP) with only binary variables, no explicit constraints, and a quadratic objective, which in turn is a subfamily of Mixed-integer nonlinear programming (MINLP), extending MIP to allow general nonlinear objective functions and constraints. By this token, QUBOs can enjoy the rich set of methods and solver capabilities that have been steadily expanding over the past three decades, from mixed-integer linear programs (MILPs), to second-order cone programs (MISOCPs), to nonconvex quadratic problems (MIQP/MIQCP), and finally to general MINLPs~\cite{bixby2012brief,koch2022progress,gurobi}, to solve these problems to global optimality.  

MINLPs are, in general, hard to solve as they often involve nonconvexities such as bilinear terms, nonlinear constraints, or logic-based disjunctions. These break convexity and frequently require reformulations, decompositions, or restrictions to tractable subclasses~\cite{kronqvist2019review}. Most approaches to solving MINLPs to global optimality rely on the branch-and-bound (B\&B) algorithm, enhanced with methods such as McCormick relaxations, the reformulation-linearization technique (RLT) cuts, and second-order cone approximations~\cite{rehfeldt2023faster}. Modern commercial solvers such as \texttt{Gurobi} incorporate these techniques and can solve QUBOs exactly to optimality. However, exact methods quickly become intractable as the problem size grows.

To handle such large instances, heuristic solvers are commonly used. Metaheuristics, such as simulated annealing, tabu search, genetic algorithms, large neighborhood search, and path relinking, can produce high-quality solutions for problems with tens or hundreds of thousands of variables. More recently, physics-inspired methods such as coherent Ising machines, digital annealers, and quantum annealers have shown promise in efficiently sampling low-energy states~\cite{mohseni2022ising}.  
Although these methods do not guarantee global optimality, they offer practical performance in large-scale QUBO instances.

This work explores hybrid approaches that pair these solvers with quantum and physics-inspired heuristics to accelerate B\&B while preserving global optimality guarantees. That is, providing both a solution and a certificate of global optimality.

\subsection*{Motivation for Hybrid Algorithms}

Although classical exact methods such as branch-and-bound (B\&B) can solve small to medium-sized QUBOs to global optimality, 
they become computationally impractical for large problems due to the exponential growth in sub-problems and the cost of solving relaxations. This motivates using quantum heuristics, thanks to the connection between QUBOs and the transverse field Ising mode~\cite{lucas2014ising}. However, current quantum hardware is limited by qubit count, sparse connectivity, and environmental noise~\cite{abbas2024challenges}. In fact, practical instances often exceed hardware capacity, which frequently prevents direct execution, or the obtained solutions are of poor quality.
More importantly, global optimality can usually not be guaranteed by quantum hardware.

These complementary limitations motivate hybrid strategies that combine classical guarantees with quantum or physics-inspired heuristics. Quantum heuristics can guide branching, supply high-quality incumbents, or tighten bounds, while the classical solver retains responsibility for global optimality certificates.  
The B\&B tree structure naturally supports this integration by enabling selective oracle calls on promising subproblems.

\subsection*{Our Contribution}

We present a practical hybrid quantum-classical algorithm for solving QUBO problems by embedding Ising-based heuristics into a classical B\&B framework.  
Our method interfaces externally with \texttt{Gurobi} and invokes classical or quantum Ising solvers as heuristic oracles in selected subproblems. These solvers are invoked at different levels of the B\&B tree.
Oracle calls may provide branching decisions or candidate incumbent solutions.
To improve effectiveness, we applied a preprocessing step to reduce the size of embedded subproblems and used a graph-aware branching rule based on variable degrees in the QUBO interaction graph.  
This design maximizes the value of Ising oracle calls while remaining within hardware constraints.  
Our implementation integrates with simulated and hardware Ising solvers and is publicly released\footnote{\url{https://github.com/SECQUOIA/QuantumBranchAndBound}}, along with a benchmark suite of over 5,000 QUBO instances from QUBOLib.jl\footnote{\url{https://github.com/SECQUOIA/QUBOLib.jl}}.

\begin{itemize}
\item \textbf{Algorithmic contributions}: (i) a decision mechanism that triggers an Ising oracle at different stages of the branch-and-bound process; (ii) a graph-aware branching rule coupled with edge-contraction preprocessing.
\item \textbf{Empirical contributions}: an open-source implementation and benchmark over more than a thousand instances show a median 17\% node reduction and 11\% wall-clock speed-up.  
The shifted geometric mean with a 10-second shift (SGM10) reduces the baseline solve time from 154 to 137 seconds.
\end{itemize}

The remainder of this paper is organized as follows.  
Section~\ref{sec:background} reviews exact, heuristic, and quantum approaches.  
Section~\ref{sec:method} describes the hybrid algorithm.  
Section~\ref{sec:results} presents empirical results.

\subsection*{Related Work}

\paragraph{Quantum-centric hybrid methods.}  
Theoretical groundwork for the integration of quantum routines into classical branch-and-bound (B\&B) was laid in~\cite{montanaro2020quantum}, which formalized the \emph{Quantum Branch and Bound} (QBB) model and proved potential speedups when quantum subroutines tighten bounds or guide branching.  
Further theoretical developments appear in~\cite{chakrabarti2022universal}.  
Building on these ideas, several prototype implementations have been proposed.  
The work in~\cite{sanavio2024hybrid} inserts D-Wave quantum-annealing calls at selected B\&B nodes to refine incumbents on small QUBO benchmarks.  
The method in~\cite{matsuyama2024efficient} proposes Quantum Relaxation-Based B\&B (QR-BnB), where a gate-based device estimates ground-state energies to tighten lower bounds.  
Reference~\cite{simen2025branch} introduces Branch-and-Bound Digitized Counterdiabatic Quantum Optimization (BB-DCQO), which branches on spins with high measurement uncertainty to focus the search.  
The architecture in~\cite{haner2024solving} sketches a forward-compatible design to offload entire subproblem relaxations to future quantum devices.  
Collectively, these works demonstrate two main QBB insertion points, \emph{bounding} and \emph{branching}, while differing in hardware platform, problem size, and optimality guarantees.

\paragraph{Classical-centric enhancements.}  
Classical improvements to B\&B focus on structural preprocessing, parameter tuning, and exploiting problem structure.  
The study in~\cite{rehfeldt2023faster} shows that exploiting sparsity, reducing symmetry, and tuning parameters in the open-source MIP solver SCIP can significantly accelerate MaxCut (a special case of QUBO) on large, sparse graphs.

\paragraph{Our approach.}  
We target currently available hardware and avoid modifying the internals of the solver.  
Rather than embedding quantum logic into the solver, we treat classical and quantum Ising engines as \emph{external oracles} invoked through callbacks.  
These oracles provide candidate incumbents and branching guidance, while the commercial solver remains responsible for bounding, node selection, and global optimality certificates.  
This black-box strategy differs from QBB approaches that require the hard-coding of quantum subroutines or the use of custom solver forks.  
It also enables a fair, large-scale empirical study in thousands of diverse QUBO instances.


\section{Background}
\label{sec:background}

This section provides a concise background on the three methodological pillars relevant to hybrid optimization: exact B\&B, classical heuristics, and quantum and physics-inspired methods.

\subsection{Exact Methods: Branch-and-Bound and Variants}

The branch-and-bound (B\&B) method is a recursive divide-and-conquer algorithm to solve mixed-integer optimization problems (MIPs) to global optimality.  
It systematically maintains upper and lower bounds on the optimal value and prunes regions of the search space that cannot contain better solutions.  
As a result, B\&B is an exact method that can find and certify optimal solutions under general assumptions.  
We refer to~\cite{conforti2014integer} for mixed-integer linear programming and to~\cite{burer2012non} for nonlinear generalizations.  
The typical steps of the B\&B algorithm are:

\begin{enumerate}
\item \textbf{Initialization (Root Node)}: Solve the continuous relaxation of the original MIP, ignoring the integrality constraints.  
This yields an initial bound on the optimal objective value.  
If the solution is already integer feasible, it is globally optimal.

\item \textbf{Branching}: If the relaxation has fractional values for integer variables, choose one such variable and create subproblems by imposing disjunctions (e.g., $x_i \leq \lfloor x_i^* \rfloor$ or $x_i \geq \lceil x_i^* \rceil$), thereby partitioning the feasible region.

\item \textbf{Bounding}: Solve the continuous relaxation for each subproblem.  
The solution provides an optimistic (upper or lower, depending on the optimization direction) bound on the objective within that subregion.

\item \textbf{Fathoming (Pruning)}: Discard a node if: (i) its relaxation is infeasible, (ii) its bound is worse than the incumbent, or (iii) its relaxed solution is integer feasible.  
If the new feasible solution improves the incumbent, update the incumbent.

\item \textbf{Node Selection}: Select the next active node for exploration.  
Strategies include best-bound, depth-first, or breadth-first.  
The goal is to explore promising parts of the tree while managing memory.

\item \textbf{Iteration}: Repeat branching, bounding, and fathoming on the selected node.  
This recursive process explores the full tree or continues until optimality is proven.

\item \textbf{Termination}: Stop when no active nodes remain.  
The incumbent is then the optimal solution.  
If no feasible solution was found, the problem is infeasible.
\end{enumerate}

While the branch-and-bound framework applies to both linear and nonlinear MIPs, solving the relaxations at each node can be computationally intensive.  
These relaxations are typically continuous and convex.  
For MILPs, they are linear programs (LPs), but more general MIPs may involve convex nonlinear relaxations.  
In some cases, these relaxations may even be nonconvex.  
Subproblems created during branching can be nearly as difficult to solve as the original problem.  
This observation has motivated the use of heuristics to either bypass expensive relaxations or produce feasible solutions that accelerate the search.  
Three common B\&B variants that exploit this idea are Branch-and-Prune, Branch-and-Cut, and Primal Heuristics~\cite{lucena1996branch}.

\begin{enumerate}
\item[\textbf{Branch and Prune:}]  
This variant focuses on recursive branching and pruning of subregions that cannot produce better solutions than the incumbent.  
Pruning decisions are based on bounds computed from continuous relaxations.  
A node is discarded if its relaxation is infeasible, yields a bound worse than the incumbent, or results in an integer-feasible solution.  
Tighter bounds enable earlier pruning and improve convergence.

\item[\textbf{Branch and Cut:}]  
This method augments B\&B by adding cutting planes to the relaxed subproblems.  
These valid inequalities are violated by the current relaxation, but they are respected by all feasible integer solutions.  
Effective cuts, such as Gomory, MIR, or cover inequalities, can tighten relaxations and improve pruning.  
Relaxations in the nodes are typically continuous and convex.  
For MILPs, this means LPs; for more general MIPs, they may be convex nonlinear programs.  
Although cutting can improve performance, excessive or poorly selected cuts may slow down per-node solution times, so cut management must be judicious.

\item[\textbf{Primal Heuristics:}]  
These methods generate high-quality feasible solutions early in the search.  
Stronger incumbent solutions improve pruning efficiency by tightening the upper bound.  
Common heuristic strategies include rounding, diving, local search, and large neighborhood search~\cite{berthold2006primal}.

\begin{itemize}
\item \textbf{Rounding}: Converts solutions from a relaxation into integer-feasible ones using rounding rules, followed by feasibility repair if needed.

\item \textbf{Diving}: Fixes variables based on their fractional values and recursively solves subproblems to explore promising regions.

\item \textbf{Local Search}: Improves a given solution by exploring its neighborhood using swap, flip, or exchange moves.

\item \textbf{Large Neighborhood Search (LNS)}: Relaxes a subset of variables and resolves the resulting reduced problem to escape local optima.

\end{itemize}
\end{enumerate}

\subsection{Heuristic Algorithms}

Given the NP-hardness of QUBO, exact methods become impractical for large-scale instances.  
To address this, the literature has proposed a wide range of heuristics and approximation algorithms capable of producing high-quality solutions in a reasonable computational time.  
In many cases, heuristics have been shown to recover optimal solutions in benchmark instances.  

\emph{Simulated Annealing (SA)} is a widely used metaheuristic inspired by the physical process of annealing in metallurgy.  
SA explores the solution space by iteratively proposing neighboring solutions~\cite{kirkpatrick1983optimization}.  
Improvements are always accepted, while worse solutions may be accepted with a probability that decreases over time according to a predefined cooling schedule.  
This mechanism helps the algorithm escape local minima and explore diverse regions of the search space.  
The performance of SA is highly sensitive to parameter tuning, especially the cooling rate.  
For QUBO problems, several effective implementations have been reported in~\cite{alkhamis1998simulated, beasley1998heuristic}.

\emph{Local Search} heuristics form another important category.  
In particular, \emph{Tabu Search} uses a short-term memory structure (the tabu list) to prevent cycling and encourage diversification.  
At each iteration, the best admissible neighbor is selected, potentially even if it worsens the objective.  
This process allows escape from local optima.  
Detailed implementations can be found in~\cite{beasley1998heuristic, glover1998adaptive}.  
\emph{Iterated Local Search} combines local search with strategic perturbations.  
Effective variants for QUBO are given in~\cite{palubeckis2004multistart}.

\emph{Genetic Algorithms (GAs)} represent a population-based heuristic that mimics evolutionary processes.  
A pool of candidate solutions evolves through crossover (recombination) and mutation (random perturbation).  
Selection is based on fitness, typically measured by the QUBO objective value.  
Many GA implementations incorporate local search to refine offspring, although convergence can be slow and quality may vary.  
These methods are computationally demanding and sensitive to population diversity.  
Relevant studies include~\cite{duarte2005low, hasan2000comparison, merz1999genetic}.

Additional heuristic families include \emph{path relinking}~\cite{festa2002randomized}, \emph{cross-entropy methods}~\cite{laguna2009hybridizing}, \emph{global equilibrium search}~\cite{pardalos2008global}, and \emph{greedy construction strategies}~\cite{festa2002randomized}.  
For a comprehensive empirical comparison across heuristics for QUBO and Max-Cut, we refer to the systematic evaluation in~\cite{dunning2018works}.

\subsection{Quantum and Physics-Inspired hybrid Methods}

Several quantum computing paradigms and their corresponding hardware implementations have been proposed to tackle difficult optimization problems. A notable example is \emph{Adiabatic Quantum Computing (AQC)}, which begins with a quantum system in the ground state of an initial Hamiltonian and slowly evolves it toward a cost Hamiltonian encoding the optimization objective.  
If the evolution is sufficiently slow, the adiabatic theorem suggests that the system remains in the ground state of the final Hamiltonian, corresponding to the optimal solution.  
However, determining how slow is ``sufficiently'' slow depends on the minimum spectral gap during evolution, which is generally intractable to compute.  
As a result, real implementations use heuristics for schedule selection and face challenges such as thermal noise, decoherence, and limited qubit connectivity.  
These limitations often require embedding logical variables using multiple physical qubits, which adds overhead.  
\emph{Quantum Annealing (QA)} models AQC in the presence of such physical imperfections and is used as a heuristic optimization method.

An alternative to AQC is the gate-based quantum computing model.  
In this framework, \emph{Variational Quantum Algorithms (VQAs)} have emerged as a family of hybrid quantum-classical methods suitable for optimization~\cite{cerezo2021variational}.  
VQAs rely on parameterized quantum circuits whose performance is evaluated by a classical optimizer based on a measured cost function.  
The classical optimizer updates the parameters to minimize this objective, typically through iterative feedback.  
Training these circuits is known to be NP-hard in general~\cite{bittel2021training}.

Among VQAs, the \emph{Quantum Approximate Optimization Algorithm (QAOA)} has become a widely studied strategy for combinatorial problems~\cite{willsch2020benchmarking}.  
QAOA alternates between applying the cost and mixing Hamiltonians.  
The number of alternations determines the circuit depth.  
Each round involves optimizing a set of continuous parameters, known as rotation angles, that control the unitary operations.  
Although QAOA is provably optimal in the infinite-depth limit due to its equivalence to AQC, its performance at finite depth remains difficult to analyze due to quantum many-body interactions and classical optimization difficulties~\cite{uvarov2021barren}.

We point out that all QAOA~\cite{farhi2014quantum}, together with Ising-based hardware such as D-Wave quantum annealers~\cite{johnson2011quantum}, and coherent Ising machines~\cite{honjo2021100} are specifically designed to tackle QUBOs heuristically. However, hardware noise and embedding overhead still limit scale, but empirical gains on medium-sized problems are encouraging~\cite{mohseni2022ising}. For a broader review of quantum heuristics for Ising problems, we refer the reader to~\cite{sanders2020compilation}.

As pointed out in the introduction, Quantum Branch-and-Bound (QBB) frameworks integrate such routines into classical B\&B, aiming to accelerate bounding or branching while preserving global optimality guarantees~\cite{matsuyama2024efficient,simen2025branch,sanavio2024hybrid}.

\section{Proposed methods}
\label{sec:method}

In this work, we propose and implement a hybrid quantum branch-and-bound (B\&B) algorithm specifically designed to solve QUBO problems. The core idea is to incorporate heuristic solutions obtained from quantum hardware into the B\&B tree to tighten the upper bound and enhance pruning efficiency. In general, modern B\&B solvers allow external solution information to be injected in three ways: MIPStart (also known as warm start), heuristic callbacks, and variable hints. Since quantum solvers typically provide complete feasible solutions rather than partial guidance, our method focuses on the first two mechanisms and does not consider variable hints, which are better suited for soft guidance rather than hard feasible inputs.

Algorithm~\ref{algo:quantum branch and cut} shows the high-level pseudocode of our method. Compared to the standard B\&B algorithm, our approach introduces three key enhancements. First, quantum solutions are injected at the root node using the MIPStart mechanism, allowing the solver to begin with a high-quality incumbent and prune large portions of the tree early on. Second, we extend this injection strategy to subtrees by invoking heuristic callbacks at interior nodes. This enables the algorithm to continually benefit from quantum-generated solutions throughout the search. Although hybrid quantum solvers can handle QUBO problems that exceed the size limits of quantum annealers, the quality and efficiency of quantum solutions tend to degrade with increasing problem size. To better exploit quantum hardware, we design branching strategies that prioritize subproblems likely to be smaller, and thus more amenable to high-quality quantum solutions. These methods are implemented in our library and extensively tested on thousands of QUBO instances. More details of the experimental results are presented in the next section.

\begin{algorithm}
\caption{Quantum Branch-and-Bound Framework}
\label{algo:quantum branch and cut}
\begin{algorithmic}[1]
\State \textbf{Inject solution} \Comment{via MIPStart}
\State \textbf{Calculate branch priority} \Comment{based on Q matrix}
\State Perform presolve
\State Solve root node LP relaxation
\While{termination criteria not met}
    \State Node selection
    \State \textbf{Inject solution}  \Comment{via heuristic callback}
    \State Node presolve
    \State Solve the LP relaxation
    \State Apply cutting planes
    \State Apply primal heuristics
    \If{a feasible integer solution is found}
        \State Update incumbent solution
    \ElsIf{the node is still feasible}
        \State Branch on fractional variables
        \State Insert child nodes into the search tree
    \EndIf
\EndWhile
\end{algorithmic}
\end{algorithm}

\subsection{Root node: MIPstart}

Quantum solvers, such as quantum annealers, are designed to solve QUBO problems by reformulating them as equivalent Ising models that can be directly mapped onto quantum hardware. For large-scale QUBO instances that exceed the capacity of physical quantum annealers, hybrid solvers combine classical and quantum resources to solve the full problem without manual decomposition. This enables the direct application of quantum solvers to the original QUBO problem. Although such solvers can return feasible solutions, their quality is not always guaranteed, especially for problems with many variables or complex landscapes \cite{pusey2008adiabatic}. However, these solutions can still provide useful upper bounds for minimization problems and can be injected into the root node of a branch-and-bound solver using \texttt{MIPStart}, which allows users to supply one or multiple feasible solutions to guide the search. This process corresponds to Step1 in Algorithm \ref{algo:quantum branch and cut}.

\subsection{Injection at internal nodes}

An extension of the root-node injection idea is to insert solutions at internal nodes, finding a high-quality solution to each branch of the tree. 
These insertions can be handled by \emph{heuristic callbacks}, a technique available in modern MIP solvers that allows users to provide a feasible solution dynamically during the tree search. These callbacks are implemented in step $7$ of Algorithm \ref{algo:quantum branch and cut}. However, these calls must be handled with care, as we now discuss. Invoking heuristic solvers at every node can be prohibitively expensive and may significantly increase the overall computation time. We propose the following strategy to mitigate this issue. Since QUBO does not have constraints, a feasible solution for one node in the tree is feasible for all its child nodes. This suggests that instead of obtaining heuristic solutions at each node, we should find a large set of high-quality solutions \emph{a priori}, i.e., before starting the branch-and-bound algorithm, and store these solutions for later use. This strategy also aligns well with the high-throughput nature of quantum solvers, which are capable of producing a large volume of feasible solutions.

\subsection{Embedding and Branching Priority}

\begin{figure}[htb]
    \centering
    \includegraphics[width=0.95\linewidth]{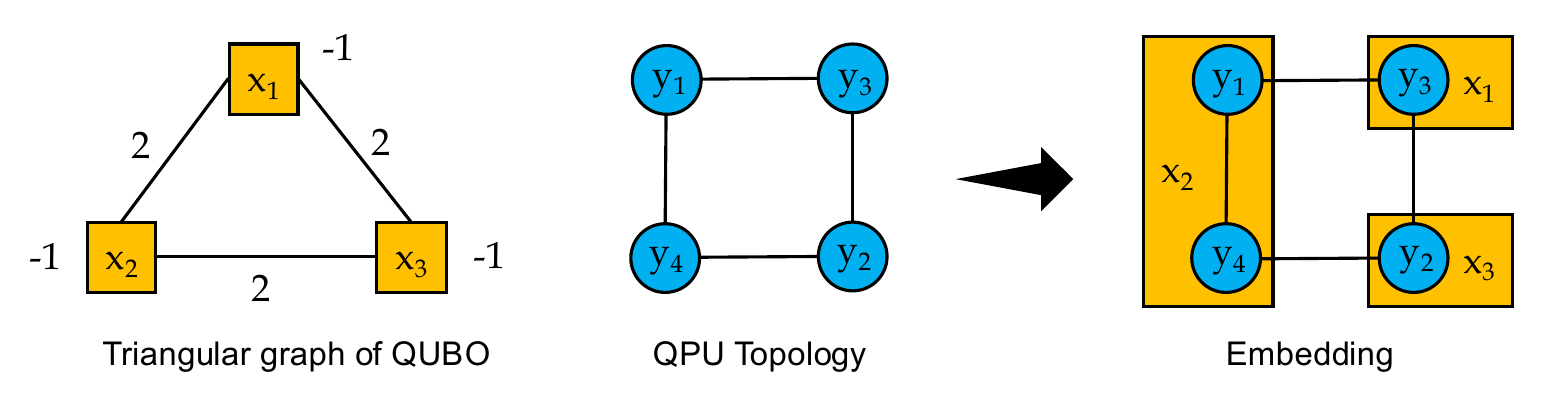}
    \caption{Example of embedding a 3-variable QUBO problem}
    \label{fig:embedding}
\end{figure}

To solve an arbitrarily posed binary quadratic problem directly on a D-Wave system requires mapping, called minor embedding, to the QPU Topology of the system's quantum processing unit (QPU) \cite{okada2019improving}. By default, D-Wave will call \texttt{minorminer} to find the embedding of the input QUBOs. Because the number of qubits inherently limits the quantum processing unit, it is desirable to embed problems with a smaller number of variables. Figure \ref{fig:embedding} shows an example of an embedding for a 3-variable QUBO \eqref{eq:embedding example} onto a four-node QPU topology \cite{D-Wave_sat_example}. The QUBO problem is first represented by the triangular graph, where nodes represent variables, and edges represent the quadratic terms. Embedding aims to map the triangular graph into the fully connected and sparse four-node graphs.

\begin{equation}
    \label{eq:embedding example}
    \min 2 x_1 x_2 + 2 x_1 x_3 + 2 x_2 x_3 - x_1 - x_2 - x_3
\end{equation}

\begin{figure}
    \centering
    \includegraphics[width=0.95\linewidth]{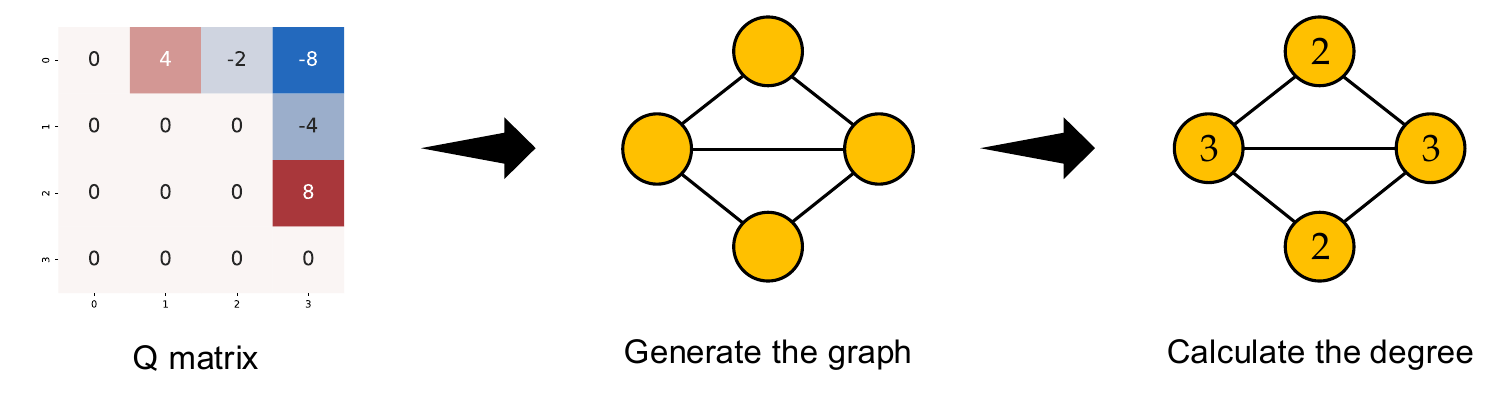}
    \caption{Example of calculating branch priority from quadratic objective matrix}
    \label{fig:priority}
\end{figure}

This suggests that one should call a heuristic solver precisely on the nodes where most variables have been fixed, which are either deep nodes (i.e., further down the tree) or nodes where variables that appear in a large number of quadratic terms have been fixed. To illustrate, consider the following example. Suppose that we are given a QUBO with a Q matrix as in Figure \ref{fig:priority}, which offers the optimization problem

\begin{equation}
\label{eq:priority example}
\min_{x \in \{0,1\}^4 } 4x_1x_2 - 2x_1x_3  - 8 x_1x_4 - 4x_2 x_4 + 8 x_3 x_4
\end{equation}

Observe now that if $x_1$ and $x_4$ are fixed, we obtain an optimization program on variables $x_2$ and $x_3$ without quadratic terms. This toy example suggests that a good rule of thumb consists of first branching on variables that appear in many quadratic terms in the objective function, as they have a larger potential to diminish the size of the problem when branching on them. 

Formally, we consider the matrix graph $G(Q)$ of $Q$ where the vertices correspond to the variables in Program \eqref{qubo}, and we add an edge between two vertices $i$ and $j$ if $Q_{ij} \neq 0$. The $degree$ of a vertex $d(i)$ is the number of edges incident to it, or equivalently, the number of quadratic terms of the form $Q_{ij}x_ix_j , j \in {N}$ in which the variable $x_i$ appears in the objective. 
We define the \emph{Branch priority} of a vertex of $G(Q)$ as its degree. Figure \ref{fig:priority} shows an example of computing branch priority for the QUBO problem \eqref{eq:priority example}. In the branching step of our proposed B\&B algorithm, we continue the iteration on the branch with the highest branch priority, breaking ties arbitrarily. This step is implemented in lines $2$ and $15$ of Algorithm \ref{algo:quantum branch and cut}.

\section{Numerical Experiments}
\label{sec:results}

\begin{figure}[htb]
    \centering
    \includegraphics[width=0.9\linewidth]{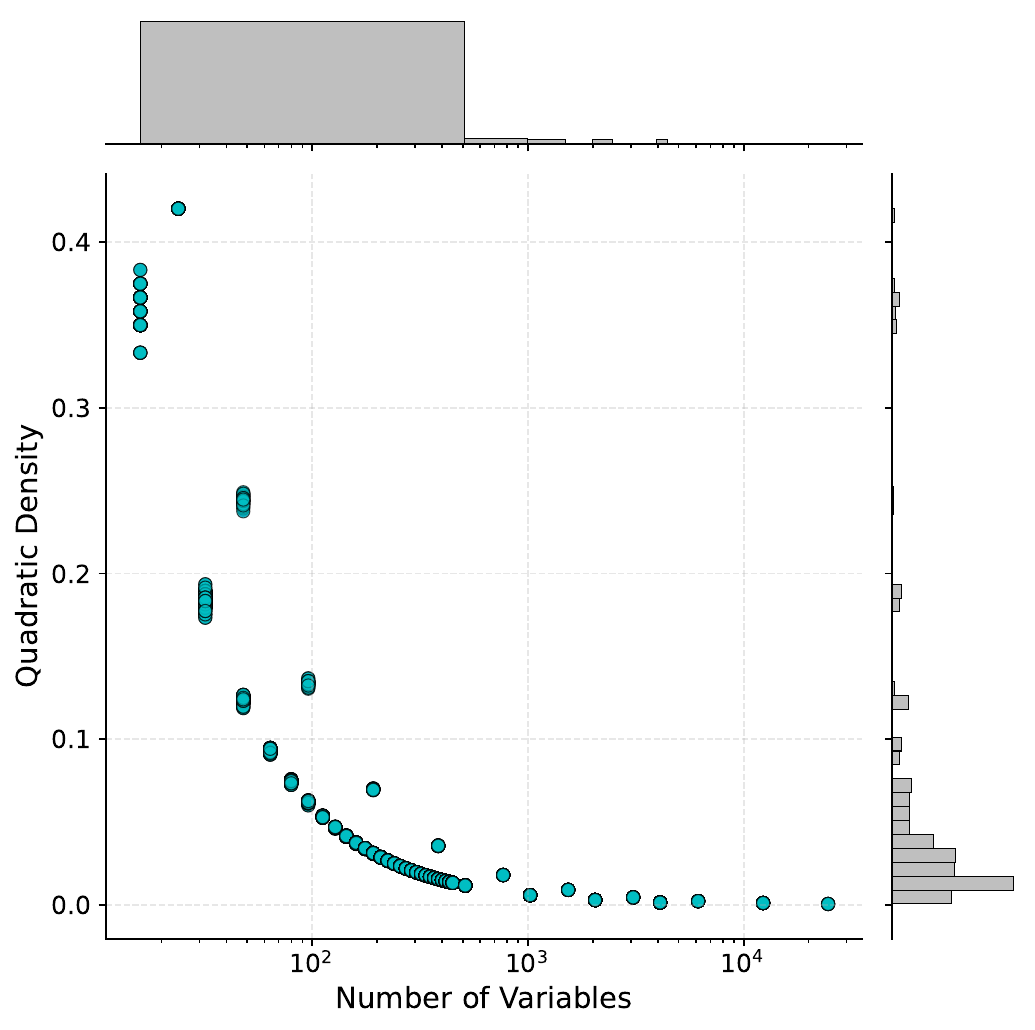}
    \caption{Quadratic sparsity of QUBOLib instances}
    \label{fig:quadratic sparsity}
\end{figure}

We evaluated the proposed quantum branch-and-bound method on a benchmark set of 5807 instances from QUBOLib, which includes planted solutions to 3-regular 3-XORSAT and 5-regular 5-XORSAT problems. Figure \ref{fig:quadratic sparsity} and Table 1 provide a statistical overview of these instances.
Figure 3 illustrates the quadratic sparsity of QUBOLib instances, revealing a clear trend: As the number of variables increases, the quadratic term sparsity also increases. 
Table 1 categorizes the entire set of benchmarks into three collections. The 3-regular 3-XORSAT problems are sourced from two different arXiv datasets, while the 5-regular 5-XORSAT problems cover significantly larger problem sizes, with up to 24,576 variables. 

Our Quantum Branch-and-Bound algorithm is implemented using a modular and extensible Julia-based pipeline. We begin by loading over 5000 QUBO instances from the QUBOLib benchmark using the \texttt{QUBOLib.jl} and \texttt{QUBOTools.jl} \cite{xavier2023qubo} packages.
Each instance is translated into a structured model using the JuMP modeling language.

To guide branching decisions within the solver, we calculate the degree of the graph induced by the quadratic objective matrix using \texttt{Graphs.jl}. The branching priority information is passed to Gurobi 11.0.0, utilizing its Branch \& Cut capabilities for exact optimization. The NonConvex parameter is set to 2 to enable solving non-convex quadratic programs, which are reformulated into bilinear forms and handled via spatial branching. All experiments are run with \texttt{ThreadLimit} = 1 and a time limit of 900 seconds.

In addition to exact methods, we incorporate heuristic warm-starts from quantum solvers such as D-Wave via MQLib \cite{dunning2018works}, allowing the solver to initialize with high-quality feasible solutions. 
We test 16 heuristic methods from MQLib , including \texttt{BURER2002}, \texttt{FESTA2002GVNSPR}, \texttt{PALUBECKIS2004bMST3}, \texttt{PALUBECKIS2006}, \texttt{FESTA2002GPR}, \texttt{FESTA2002GVNS}, \texttt{MERZ2004}, \texttt{PALUBECKIS2004bMST2}, \texttt{BEASLEY1998TS}, \texttt{LU2010}, \texttt{FESTA2002G}, \texttt{PALUBECKIS2004bMST1}, \texttt{MERZ1999GLS}, \texttt{MERZ2002KOPT}, \texttt{ALKHAMIS1998}, \texttt{MERZ2002GREEDYKOPT}.
We also evaluate simulated annealing using \texttt{dimod.neal} (v0.5.9), and quantum annealing on D-Wave's Advantage 4.1 system with 5,750 qubits and over 35,000 couplers.

To focus our analysis, we filter the dataset to instances that (i) take more than 10 seconds to solve using default Gurobi, and (ii) can be solved to optimality by at least one of the tested methods within the time limit. This yields a refined test set of 1,454 instances for detailed comparison and analysis.

We measure performance using the shifted geometric mean (SGM) of solve time and number of explored nodes, with a shift of 10 (SGM10). 
If the instance is not solved to optimality within the time limit, the solve time is always set to the corresponding time limit, and we record the number of explored nodes. The results are presented in Table \ref{tab:result_summary}.

In particular, MQLib returns only the best-found solution to the given problem. When used with Gurobi's MIP start strategy, this solution is injected at the root node. Alternatively, when used in a heuristic callback strategy, MQLib is invoked at every node to attempt to solve subproblems. For SA and QA, we experiment with injecting the top 1, 10, 30, or 100 solutions sorted by objective value. In callback mode, SA and QA are applied once at the root node to generate a solution pool. During the branch-and-bound process, solutions from this pool are selectively injected based on the node subproblem. Moreover, to test the upper bound of the improvement, we tested the performance of providing the best solution in the MIP start strategy.

The results are summarized in Table \ref{tab:result_summary}. Among the 16 MQLib heuristics tested, \texttt{PALUBECKIS2006} consistently achieves the best performance. To simplify the presentation, we only report the results of \texttt{PALUBECKIS2006} in Table \ref{tab:result_summary} as the representative MQLib method.
It is shown that using branch priority alone improves Gurobi's performance by 17.3\% in node count and 11.1\% in runtime. When using MIP, start with \texttt{PALUBECKIS2006}, simulated annealing, or quantum annealing, we observe approximately a 10\% runtime improvement. For simulated annealing, injecting more solutions leads to a modest 3\% additional improvement, while quantum annealing shows limited sensitivity to the number of solutions provided. Combining MIP start with branch priority yields better results than MIP start alone, but still slightly underperforms the branch priority strategy alone. 
The results of injecting the best solution demonstrate the upper bound of improvement achievable via solution injection: the runtime is reduced by 83.0\%, and the number of explored nodes drops by 90.6\%, while still solving 1170 out of 1454 instances. Although the results remain similar when combined with priority, they confirm the potential of how high-quality starts can dramatically accelerate problem solving.

When applying heuristic callbacks, we find that invoking MQLib at every node introduces significant overhead, leading to longer runtimes and a reduced solve rate of 921 out of 1,454 instances. Although SA and QA callbacks are applied more efficiently and invoked only once, they still result in longer solve times and slightly fewer explored nodes. These findings suggest that node-level heuristic injection is often too costly in practice and should be used with caution.

\begin{table}[htbp]
\scriptsize
\centering
\label{tab:qubolib summary}
\begin{tabular}{cccc}
\toprule
\textbf{Collection} & \textbf{\# of instances} & \textbf{\# of variables} \\
\midrule
3-Regular 3-XORSAT \cite{kowalsky20223} & 2300 & 16 $\sim$ 4096 \\
\addlinespace
3-Regular 3-XORSAT \cite{hen2019equation} & 3200 & 16 $\sim$ 4096 \\
\addlinespace
5-Regular 5-XORSAT \cite{hen2019equation} & 307 & 24 $\sim$ 24576 \\
\bottomrule
\end{tabular}
\caption{XORSAT Planted Solutions Collections - QUBOLib}
\end{table}

\begin{table*}[htbp]
\centering
\label{tab:result_summary}
\resizebox{\textwidth}{!}{
\begin{tabular}{llrll}
\toprule
Strategy & Heuristic & \# solved instances & Node Count  & Runtime [s] \\
\midrule
Baseline & - & 1065 & 28170.0 & 154.5 \\
 Branch Priority  & - & 1077 & 23304.8 (-17.3\%) & 137.4 (-11.1\%) \\
 MIP Start  & MQLib PALUBECKIS2006 & 1073 & 24543.0 (-12.9\%) & 142.5 (-7.8\%) \\
 MIP Start  & DWave SA TOP1 & 1070 & 25418.8 (-9.8\%) & 150.1 (-2.9\%) \\
 MIP Start  & DWave SA TOP10 & 1080 & 25256.3 (-10.3\%) & 148.3 (-4.1\%) \\
 MIP Start  & DWave SA TOP30 & 1074 & 24616.1 (-12.6\%) & 144.8 (-6.3\%) \\
 MIP Start  & DWave SA TOP100 & 1059 & 24665.1 (-12.4\%) & 145.1 (-6.1\%) \\
 MIP Start  & DWave QA TOP1 & 1081 & 24958.7 (-11.4\%) & 148.8 (-3.7\%) \\
 MIP Start  & DWave QA TOP10 & 1074 & 24880.3 (-11.7\%) & 149.8 (-3.1\%) \\
 MIP Start  & DWave QA TOP30 & 1069 & 25069.6 (-11.0\%) & 150.1 (-2.8\%) \\
 MIP Start  & DWave QA TOP100 & 1079 & 25189.0 (-10.6\%) & 150.6 (-2.6\%) \\
 MIP Start + Embedding  & DWave QA TOP1 & 1073 & 25713.4 (-8.7\%) & 148.2 (-4.1\%) \\
 MIP Start + Embedding  & DWave QA TOP10 & 1078 & 24742.1 (-12.2\%) & 142.9 (-7.5\%) \\
 MIP Start + Embedding  & DWave QA TOP30 & 1077 & 25093.5 (-10.9\%) & 146.0 (-5.5\%) \\
 MIP Start + Embedding  & DWave QA TOP100 & 1069 & 25536.7 (-9.3\%) & 146.6 (-5.1\%) \\
 MIP Start  & Best Solution & 1170 & 2643.3 (-90.6\%) & 26.3 (-83.0\%) \\
  MIP Start + Branch Priority  & MQLib PALUBECKIS2006 & 1094 & 23292.9 (-17.3\%) & 137.7 (-10.9\%) \\
 MIP Start + Branch Priority  & DWave SA TOP1 & 1068 & 23440.7 (-16.8\%) & 138.0 (-10.7\%) \\
 MIP Start + Branch Priority  & DWave SA TOP10 & 1079 & 23956.5 (-15.0\%) & 140.3 (-9.2\%) \\
 MIP Start + Branch Priority  & DWave SA TOP30 & 1083 & 24000.7 (-14.8\%) & 140.5 (-9.1\%) \\
 MIP Start + Branch Priority  & DWave SA TOP100 & 1077 & 23550.5 (-16.4\%) & 138.1 (-10.6\%) \\
 MIP Start + Branch Priority  & DWave QA TOP1 & 1077 & 24312.4 (-13.7\%) & 142.3 (-7.9\%) \\
 MIP Start + Branch Priority  & DWave QA TOP10 & 1072 & 23628.1 (-16.1\%) & 139.0 (-10.0\%) \\
 MIP Start + Branch Priority  & DWave QA TOP30 & 1093 & 23714.2 (-15.8\%) & 138.9 (-10.1\%) \\
 MIP Start + Branch Priority  & DWave QA TOP100 & 1071 & 24342.5 (-13.6\%) & 141.4 (-8.5\%) \\
 MIP Start + Branch Priority  & Best Solution & 1164 & 2649.6 (-90.6\%) & 26.5 (-82.9\%) \\
 Heuristic Callback  & MQLib PALUBECKIS2006 & 921 & 17698.5 (-37.2\%) & 307.8 (+99.2\%) \\
 Heuristic Callback  & DWave QA  & 1060 & 27785.0 (-1.4\%) & 160.5 (+3.9\%) \\
 Heuristic Callback  & DWave SA  & 1039 & 26559.0 (-5.7\%) & 183.9 (+19.0\%) \\
\bottomrule
\end{tabular}}
\caption{Summary of solver performance on 1454 instances in QUBOLib (SGM10)}
\end{table*}

\section{Conclusions}

This work proposes a practical hybrid quantum-classical branch-and-bound framework for solving QUBO problems to global optimality. The proposed method provides a unified framework to integrate Ising solvers, including both classical heuristics and quantum annealers, into a Gurobi-based branch-and-bound solver. Extensive experiments on over 5,800 instances from QUBOLib show that warm-starting with high-quality solutions from Ising solvers yields a 5\% improvement, and a carefully designed branch priority rule alone can reduce solve time and node count by over 10\%. However, the improvement remains well below the potential upper bound obtained by providing the best solution. Additionally, node-wise heuristic callbacks are computationally expensive and often counterproductive. Overall, our results validate the potential of hybrid quantum-classical strategies to accelerate exact solvers on structured QUBO problems. Developing more effective methods for integrating quantum solvers as node-wise heuristics remains an open direction for future research.

\section{ACKNOWLEDGMENTS}
D.B.N. and Z.P. acknowledge the support of the startup grant from the Davidson School of Chemical Engineering at Purdue University.
This material is based upon work supported by the Center for Quantum Technologies under the Industry-University Cooperative Research Center Program at the US National Science Foundation under Grant No. 2224960.

\bibliographystyle{alpha}
\bibliography{references}

\newcommand{\etalchar}[1]{$^{#1}$}
\begin{thebibliography}{WWJ{\etalchar{+}}20}

\bibitem[AAA{\etalchar{+}}24]{abbas2024challenges}
Amira Abbas, Andris Ambainis, Brandon Augustino, Andreas B{\"a}rtschi, Harry Buhrman, Carleton Coffrin, Giorgio Cortiana, Vedran Dunjko, Daniel~J Egger, Bruce~G Elmegreen, et~al.
\newblock {Challenges and opportunities in quantum optimization}.
\newblock {\em Nature Reviews Physics}, pages 1--18, 2024.

\bibitem[AHA98]{alkhamis1998simulated}
Talal~M Alkhamis, Merza Hasan, and Mohamed~A Ahmed.
\newblock {Simulated annealing for the unconstrained quadratic pseudo-Boolean function}.
\newblock {\em European journal of operational research}, 108(3):641--652, 1998.

\bibitem[Bea98]{beasley1998heuristic}
John~E Beasley.
\newblock {Heuristic algorithms for the unconstrained binary quadratic programming problem}, 1998.

\bibitem[Ber06]{berthold2006primal}
Timo Berthold.
\newblock {\em {Primal heuristics for mixed integer programs}}.
\newblock PhD thesis, Zuse Institute Berlin (ZIB), 2006.

\bibitem[Bix12]{bixby2012brief}
Robert~E Bixby.
\newblock {A brief history of linear and mixed-integer programming computation}.
\newblock {\em Documenta Mathematica}, 2012:107--121, 2012.

\bibitem[BK21]{bittel2021training}
Lennart Bittel and Martin Kliesch.
\newblock {Training variational quantum algorithms is NP-hard}.
\newblock {\em Physical review letters}, 127(12):120502, 2021.

\bibitem[BL12]{burer2012non}
Samuel Burer and Adam~N Letchford.
\newblock {Non-convex mixed-integer nonlinear programming: A survey}.
\newblock {\em Surveys in Operations Research and Management Science}, 17(2):97--106, 2012.

\bibitem[CAB{\etalchar{+}}21]{cerezo2021variational}
Marco Cerezo, Andrew Arrasmith, Ryan Babbush, Simon~C Benjamin, Suguru Endo, Keisuke Fujii, Jarrod~R McClean, Kosuke Mitarai, Xiao Yuan, Lukasz Cincio, et~al.
\newblock {Variational quantum algorithms}.
\newblock {\em Nature Reviews Physics}, 3(9):625--644, 2021.

\bibitem[CCZ{\etalchar{+}}14]{conforti2014integer}
Michele Conforti, G{\'e}rard Cornu{\'e}jols, Giacomo Zambelli, Michele Conforti, G{\'e}rard Cornu{\'e}jols, and Giacomo Zambelli.
\newblock {\em {Integer programming models}}.
\newblock Springer, 2014.

\bibitem[CMYP22]{chakrabarti2022universal}
Shouvanik Chakrabarti, Pierre Minssen, Romina Yalovetzky, and Marco Pistoia.
\newblock {Universal quantum speedup for branch-and-bound, branch-and-cut, and tree-search algorithms}.
\newblock {\em arXiv preprint arXiv:2210.03210}, 2022.

\bibitem[DGS18]{dunning2018works}
Iain Dunning, Swati Gupta, and John Silberholz.
\newblock {What works best when? A systematic evaluation of heuristics for Max-Cut and QUBO}.
\newblock {\em INFORMS Journal on Computing}, 30(3):608--624, 2018.

\bibitem[DSFC05]{duarte2005low}
Abraham Duarte, {\'A}ngel S{\'a}nchez, Felipe Fern{\'a}ndez, and Ra{\'u}l Cabido.
\newblock {A low-level hybridization between memetic algorithm and VNS for the max-cut problem}.
\newblock In {\em Proceedings of the 7th annual conference on Genetic and evolutionary computation}, pages 999--1006, 2005.

\bibitem[FGG14]{farhi2014quantum}
Edward Farhi, Jeffrey Goldstone, and Sam Gutmann.
\newblock {A quantum approximate optimization algorithm}.
\newblock {\em arXiv preprint arXiv:1411.4028}, 2014.

\bibitem[FPRR02]{festa2002randomized}
Paola Festa, Panos~M Pardalos, Mauricio~GC Resende, and Celso~C Ribeiro.
\newblock {Randomized heuristics for the MAX-CUT problem}.
\newblock {\em Optimization methods and software}, 17(6):1033--1058, 2002.

\bibitem[GKA98]{glover1998adaptive}
Fred Glover, Gary~A Kochenberger, and Bahram Alidaee.
\newblock {Adaptive memory tabu search for binary quadratic programs}.
\newblock {\em Management Science}, 44(3):336--345, 1998.

\bibitem[GKHD22]{QUBOHandbook}
Fred Glover, Gary Kochenberger, Rick Hennig, and Yu~Du.
\newblock {Quantum bridge analytics I: a tutorial on formulating and using QUBO models}.
\newblock {\em Annals of Operations Research}, 314(1):141--183, 2022.

\bibitem[{Gur}24]{gurobi}
{Gurobi Optimization, LLC}.
\newblock {Gurobi Optimizer Reference Manual}, 2024.

\bibitem[HAA00]{hasan2000comparison}
Merza Hasan, Talal Alkhamis, and Jafar Ali.
\newblock {A comparison between simulated annealing, genetic algorithm and tabu search methods for the unconstrained quadratic Pseudo-Boolean function}.
\newblock {\em Computers \& industrial engineering}, 38(3):323--340, 2000.

\bibitem[HBBZ24]{haner2024solving}
Thomas H{\"a}ner, Kyle~EC Booth, Sima~E Borujeni, and Elton~Yechao Zhu.
\newblock {Solving QUBOs with a quantum-amenable branch and bound method}.
\newblock {\em arXiv preprint arXiv:2407.20185}, 2024.

\bibitem[Hen19]{hen2019equation}
Itay Hen.
\newblock Equation planting: a tool for benchmarking ising machines.
\newblock {\em Physical Review Applied}, 12(1):011003, 2019.

\bibitem[HSI{\etalchar{+}}21]{honjo2021100}
Toshimori Honjo, Tomohiro Sonobe, Kensuke Inaba, Takahiro Inagaki, Takuya Ikuta, Yasuhiro Yamada, Takushi Kazama, Koji Enbutsu, Takeshi Umeki, Ryoichi Kasahara, et~al.
\newblock {100,000-spin coherent Ising machine}.
\newblock {\em Science advances}, 7(40):eabh0952, 2021.

\bibitem[JAG{\etalchar{+}}11]{johnson2011quantum}
Mark~W Johnson, Mohammad~HS Amin, Suzanne Gildert, Trevor Lanting, Firas Hamze, Neil Dickson, Richard Harris, Andrew~J Berkley, Jan Johansson, Paul Bunyk, et~al.
\newblock {Quantum annealing with manufactured spins}.
\newblock {\em Nature}, 473(7346):194--198, 2011.

\bibitem[KAHL22]{kowalsky20223}
Matthew Kowalsky, Tameem Albash, Itay Hen, and Daniel~A Lidar.
\newblock 3-regular three-xorsat planted solutions benchmark of classical and quantum heuristic optimizers.
\newblock {\em Quantum Science and Technology}, 7(2):025008, 2022.

\bibitem[KBLG19]{kronqvist2019review}
Jan Kronqvist, David~E Bernal, Andreas Lundell, and Ignacio~E Grossmann.
\newblock {A review and comparison of solvers for convex MINLP}.
\newblock {\em Optimization and Engineering}, 20:397--455, 2019.

\bibitem[KBPV22]{koch2022progress}
Thorsten Koch, Timo Berthold, Jaap Pedersen, and Charlie Vanaret.
\newblock {Progress in mathematical programming solvers from 2001 to 2020}.
\newblock {\em EURO Journal on Computational Optimization}, 10:100031, 2022.

\bibitem[KGJV83]{kirkpatrick1983optimization}
Scott Kirkpatrick, C~Daniel Gelatt~Jr, and Mario~P Vecchi.
\newblock {Optimization by simulated annealing}.
\newblock {\em science}, 220(4598):671--680, 1983.

\bibitem[LB96]{lucena1996branch}
Abilio Lucena and John~E Beasley.
\newblock {Branch and cut algorithms}.
\newblock {\em Advances in linear and integer programming}, 4:187--221, 1996.

\bibitem[LDM09]{laguna2009hybridizing}
Manuel Laguna, Abraham Duarte, and Rafael Marti.
\newblock {Hybridizing the cross-entropy method: An application to the max-cut problem}.
\newblock {\em Computers \& Operations Research}, 36(2):487--498, 2009.

\bibitem[Luc14]{lucas2014ising}
Andrew Lucas.
\newblock {Ising formulations of many NP problems}.
\newblock {\em Frontiers in physics}, 2:74887, 2014.

\bibitem[MF99]{merz1999genetic}
Peter Merz and Bernd Freisleben.
\newblock {Genetic algorithms for binary quadratic programming}.
\newblock In {\em Proceedings of the genetic and evolutionary computation conference}, volume~1, pages 417--424. Morgan Kaufmann Orlando, FL, 1999.

\bibitem[MHNY24]{matsuyama2024efficient}
Hiromichi Matsuyama, Wei-hao Huang, Kohji Nishimura, and Yu~Yamashiro.
\newblock {Efficient Internal Strategies in Quantum Relaxation based Branch-and-Bound}.
\newblock {\em arXiv preprint arXiv:2405.00935}, 2024.

\bibitem[MMB22]{mohseni2022ising}
Naeimeh Mohseni, Peter~L McMahon, and Tim Byrnes.
\newblock {Ising machines as hardware solvers of combinatorial optimization problems}.
\newblock {\em Nature Reviews Physics}, 4(6):363--379, 2022.

\bibitem[Mon20]{montanaro2020quantum}
Ashley Montanaro.
\newblock {Quantum speedup of branch-and-bound algorithms}.
\newblock {\em Physical Review Research}, 2(1):013056, 2020.

\bibitem[OOTT19]{okada2019improving}
Shuntaro Okada, Masayuki Ohzeki, Masayoshi Terabe, and Shinichiro Taguchi.
\newblock {Improving solutions by embedding larger subproblems in a D-Wave quantum annealer}.
\newblock {\em Scientific reports}, 9(1):2098, 2019.

\bibitem[Pal04]{palubeckis2004multistart}
Gintaras Palubeckis.
\newblock {Multistart tabu search strategies for the unconstrained binary quadratic optimization problem}.
\newblock {\em Annals of Operations Research}, 131:259--282, 2004.

\bibitem[PND08]{pusey2008adiabatic}
L~Pusey-Nazzaro and P~Date.
\newblock {Adiabatic quantum optimization fails to solve the knapsack problem. arXiv 2020}.
\newblock {\em arXiv preprint arXiv:2008.07456}, 2008.

\bibitem[PPSS08]{pardalos2008global}
Panos~M Pardalos, Oleg~A Prokopyev, Oleg~V Shylo, and Vladimir~P Shylo.
\newblock {Global equilibrium search applied to the unconstrained binary quadratic optimization problem}.
\newblock {\em Optimisation Methods and Software}, 23(1):129--140, 2008.

\bibitem[RKS23]{rehfeldt2023faster}
Daniel Rehfeldt, Thorsten Koch, and Yuji Shinano.
\newblock {Faster exact solution of sparse MaxCut and QUBO problems}.
\newblock {\em Mathematical Programming Computation}, 15(3):445--470, 2023.

\bibitem[SBC{\etalchar{+}}20]{sanders2020compilation}
Yuval~R Sanders, Dominic~W Berry, Pedro~CS Costa, Louis~W Tessler, Nathan Wiebe, Craig Gidney, Hartmut Neven, and Ryan Babbush.
\newblock {Compilation of fault-tolerant quantum heuristics for combinatorial optimization}.
\newblock {\em PRX quantum}, 1(2):020312, 2020.

\bibitem[SRC{\etalchar{+}}25]{simen2025branch}
Anton Simen, Sebasti{\'a}n~V Romero, Alejandro~Gomez Cadavid, Enrique Solano, and Narendra~N Hegade.
\newblock {Branch-and-bound digitized counterdiabatic quantum optimization}.
\newblock {\em arXiv preprint arXiv:2504.15367}, 2025.

\bibitem[STE24]{sanavio2024hybrid}
Claudio Sanavio, Edoardo Tignone, and Elisa Ercolessi.
\newblock {Hybrid classical--quantum branch-and-bound algorithm for solving integer linear problems}.
\newblock {\em Entropy}, 26(4):345, 2024.

\bibitem[UB21]{uvarov2021barren}
AV~Uvarov and Jacob~D Biamonte.
\newblock {On barren plateaus and cost function locality in variational quantum algorithms}.
\newblock {\em Journal of Physics A: Mathematical and Theoretical}, 54(24):245301, 2021.

\bibitem[WWJ{\etalchar{+}}20]{willsch2020benchmarking}
Madita Willsch, Dennis Willsch, Fengping Jin, Hans De~Raedt, and Kristel Michielsen.
\newblock {Benchmarking the quantum approximate optimization algorithm}.
\newblock {\em Quantum Information Processing}, 19:1--24, 2020.

\bibitem[XRA{\etalchar{+}}23]{xavier2023qubo}
Pedro~Maciel Xavier, Pedro Ripper, Tiago Andrade, Joaquim~Dias Garcia, Nelson Maculan, and David E~Bernal Neira.
\newblock Qubo. jl: A julia ecosystem for quadratic unconstrained binary optimization.
\newblock {\em arXiv preprint arXiv:2307.02577}, 2023.

\end{thebibliography}

\end{document}